# Arbitrarily Substantial Number Representation for Complex Number

Satrya Fajri Pratama, Azah Kamilah Muda, Yun-Huoy Choo
*Computational Intelligence and Technologies (CIT) Research Group,*
*Center of Advanced Computing and Technologies,*
*Faculty of Information and Communication Technology,*
*Universiti Teknikal Malaysia Melaka*
*Hang Tuah Jaya, 76100 Durian Tunggal, Melaka, Malaysia*
*azah@utem.edu.my*

***Abstract*— Researchers are often perplexed when their machine learning algorithms are required to deal with complex number. Various strategies are commonly employed to project complex number into real number, although it is frequently sacrificing the information contained in the complex number. This paper proposes a new method and four techniques to represent complex number as real number, without having to sacrifice the information contained. The proposed techniques are also capable of retrieving the original complex number from the representing real number, with little to none of information loss. The promising applicability of the proposed techniques has been demonstrated and worth to receive further exploration in representing the complex number.**

## I. Introduction

Machine learning (ML) process is heavily depended on the properties of its dataset. A well-constructed dataset oftentimes leads to the satisfactory performance of the ML algorithm, and vice versa. The values of the dataset attributes can be either numerical values or nominal values [1], with the numerical values as the most prominently used attributes.

For most of the times, the numerical value is a single real number. However, in some rare cases, the numerical value is a complex number. Researchers are often baffled when encountering these complex numbers. This is because ML algorithms and toolkits are commonly designed to handle real numerical value, such as WEKA [2].

There exist some strategies employed in handling the complex number. Some researchers chose to ignore the imaginary number part of the complex number and only deals with the real number part [3, 4], whereas others decide to split the complex number into 2 numerical values, effectively doubling the number of attributes used by ML algorithms [5, 6]. Some other researchers represent the complex number as a numerical or nominal value, providing a simple mapping between one value to another [7], while others decide to leave the complex number as is, and treat it as a nominal or textual data. There are also other strategies not covered in this paper.

Hence, it is evident that there are no universal or standard mechanisms on handling the complex number. Each strategy presents their own weaknesses. Ignoring the imaginary part is commendably altering the nature of the data, since there are losses of information associated with the imaginary number part of the complex number [8].

On the other hand, splitting the complex number into two values theoretically retain the information, since both numbers is intact; however, the relationship between these two values can be lost, since some ML algorithms doesn't maintain the dependencies between two attributes, most notably is feature selection algorithms [1, 9]. Feature selection algorithms may deem an attribute which contains real number part is important, and thus discarding the attribute containing imaginary number part, or vice versa.

Conversely, mapping a complex number with another numerical value, preferably natural number, or even with nominal value, present another set of challenges. The mapping will be increasingly larger when the complex number is continuous. Lastly, even though treating the complex number as a textual value seems to be the safest option, since there is no loss of information occurred, the ML algorithms may not be able to determine the similarity of one attribute to another.

Hence, it is necessary to formulate a procedure to represent the complex number in the ML domain without having to sacrifice the information contained, and thus reducing the space required to retain the information and allowing better inference to be obtained. Therefore, this paper proposes novel complex number representation algorithms which can retain the information, and more importantly, allows for original value reconstruction. This concept is almost like the feature projection method of dimensionality reduction. The remainder of the paper is structured as follows. In next section, the proposed techniques are provided. In Section 3, experimental setup involving the dataset preparation and experimental design are presented. The outcomes showcasing the reconstruction capability of the proposed techniques, and conclusion and future works are elaborated in Sections 5 and 6, respectively.

## II. Proposed Complex Number Representation

Every complex number can be expressed by specifying either the Cartesian coordinates (CC) or the polar coordinates (PC). The complex number $c$ can be represented in CC as

$$c = x + y\hat{\imath} \tag{1}$$

where $\hat{\imath}$ is the imaginary unit. The CC $x, y$ can be described as PC $r, \varphi$ with $r \geq 0$ and $\varphi \in [0, 2\pi)$ using

$$r = \sqrt{x^2 + y^2}$$
$$\varphi = \begin{cases} \operatorname{atan2}(y, x) & \operatorname{atan2}(y, x) \geq 0 \\ \operatorname{atan2}(y, x) + 2\pi & \operatorname{atan2}(y, x) < 0 \end{cases} \tag{2}$$

where $\operatorname{atan2}(y, x)$ is a special case of the arctangent function such that

$$\operatorname{atan2}(y, x) = \begin{cases} \arctan\left(\frac{y}{x}\right) & x > 0 \\ \arctan\left(\frac{y}{x}\right) + \pi & x < 0, y \geq 0 \\ \arctan\left(\frac{y}{x}\right) - \pi & x < 0, y < 0 \\ \frac{\pi}{2} & x = 0, y > 0 \\ -\frac{\pi}{2} & x = 0, y < 0 \\ \text{undefined} & x = 0, y = 0 \end{cases} \tag{3}$$

and conversely using

$$\begin{aligned} x &= r \cos \varphi \\ y &= r \sin \varphi \end{aligned} \tag{4}$$

Both ways of representing the complex numbers as a CC or as PC will resulting in twice of the features vector space and the correlation between $x, y$ or $r, \varphi$ can be lost in the ML process [8]. Thus, both two value pairs should be distinctly encoded into a single unique value, hence the correlation of the two values can be conserved. Pairing function (PF) can be employed to achieve this goal. In this paper, two renowned PFs are used, which are Cantor [10] and Szudzik [11] PFs. The formula to calculate Cantor PF and its inverse function is defined in Equations (5) and (6) respectively as

$$C = \frac{(p + q)(p + q + 1)}{2} + q \tag{5}$$

$$\begin{aligned} p &= w - q \\ q &= C - t \\ t &= \frac{w^2}{2} + 1 \\ w &= \left\lfloor \frac{\sqrt{8C} - 1}{2} \right\rfloor \end{aligned} \tag{6}$$

while the formula to calculate Szudzik PF and its inverse function is defined in Equations (7) and (8) respectively as

$$S = \begin{cases} q^2 + p & p \neq \max(p, q) \\ p^2 + p + q & p = \max(p, q) \end{cases} \tag{7}$$

$$\begin{aligned} \langle p, q \rangle &= \begin{cases} \langle S - k^2, k \rangle & S - k^2 < k \\ \langle k, S - k^2 - k \rangle & S - k^2 \geq k \end{cases} \\ k &= \lfloor \sqrt{S} \rfloor \end{aligned} \tag{8}$$

where $p$ is the first number, $q$ is the second number, $C$ is the Cantor paired value, and $S$ is the Szudzik paired value.

However, since these PFs can only be employed to distinctly encode positive natural numbers [11, 12], both $x, y$ or $r, \varphi$ which are usually kept as 64-bit double-precision floating-number format, should be transformed to natural numbers to be paired. Ref. [13] defines IEEE 754 standard to transform the double-precision floating-number format as a binary string, which consequently alterable as a long integer number. The value of a double-precision floating-number is given as

$$n = (-1)^s \left(1 + \sum_{k=1}^{52} b_{52-k} 2^{-k}\right) 2^{e-1023} \tag{9}$$

where $s$ is the sign of the floating-number, $k$ is index of bit in the binary string, and $b_i$ is the bit in the specified index. These long integer numbers then can be used as the input of PFs. $x$ and $y$ values are more favorable compared to $r$ and $\varphi$, because there is no precision lost when calculating $r$ and $\varphi$ from $x$ and $y$. However, Cantor and Szudzik PFs may only be used if both $x$ and $y$ is positive. And thus, these representations are labelled as polar Cantor and polar Szudzik representation.

To overcome these limitations, since the $x$ and $y$ or $r$ and $\varphi$ are represented in binary string, a bit-interleaved PF can also be employed. The idea is to interleave one bit from first number, followed by one bit from second number, and then followed by the subsequent bit from first number, and so on. However, both binary strings must have equal number of bits, generally achieved by concatenating zeroes in front of shorter binary string. Bit-interleaved PF has more advantages compared to Cantor and Szudzik PFs, namely the original value can be retained without any precision loss, secondly, it is applicable for any values of both $x$ and $y$ pair and $r$ and $\varphi$ pair, and lastly, it is less computationally expensive. Since it is applicable to both CC and PC, it is then named Cartesian bit-interleaved and polar bit-interleaved, respectively.

All four proposed PFs will undoubtedly produce a rather arbitrarily large representative number, at least in the range of $10^{37}$. These large numbers sometimes affect the performance of the ML algorithm, such as support vector machine and multilayer perceptron, and thus it is preferable to have a set of numbers in smaller range. Hence, the representations can be normalized in certain cases like this, preferably normalized by the value of $10^{37}$. Furthermore, it should be noted that these numbers are merely for representation purpose only, and it is not intended for arithmetic operations.

## III. Experimental Setup

In this section, a description of the experimental method is provided to conduct an extensive and rigorous study, which are the process of generating the dataset used in this study and the proposed technique validation method.

All the proposed techniques are implemented using Java 8 programming language. A total of 4,294,967,295 complex numbers are randomly generated in the range of [–9,223,372,036,854,775,807; 9,223,372,036,854,775,807). The total number is selected because it is the total amount of 32-bit integer number in Java programming language, while the range is selected because it is the minimum and maximum values of 64-bit long integer number in Java programming language.

To justify the quality of each proposed techniques, the representation values will be calculated, and from the representation values, the original complex number is retrieved. The representation quality is defined as the error between original complex number and calculated complex number from representation value. The best representation value is determined as the one with smallest error. The formula of calculating the error is defined as

$$\varepsilon = \frac{\sqrt{(x_1 - x_2)^2 + (y_1 - y_2)^2}}{r_1} \tag{10}$$

where $x_1$ and $y_1$ is the real and imaginary part of original complex number, $x_2$ and $y_2$ is the real and imaginary part of calculated complex number, and $r_1$ is the modulus of the original complex number calculated using (2).

## IV. EXPERIMENTAL RESULTS AND DISCUSSION

The representation quality of the proposed techniques is the primary consideration of this paper. As discussed earlier, a set of random complex numbers will be generated, and its corresponding representation number using the proposed techniques will be calculated. From each of the calculated representation number, the original complex number will be determined. The sample of the error calculation using (10) for a random complex number and its representations are shown in Table 1. Meanwhile, Table 2 shows the results of maximum and average representation errors of each technique.

Table 1
Sample of error calculation for a random complex number and its representations

| Source | Attribute | Value |
|---|---|---|
| Original number | Real number part | 6.7771673222051697E18 |
| | Imaginary number part | 3.6003875414142131E18i |
| | Modulus of PC | 7.67416362618991E18 |
| | Theta of PC | 0.4883359535588942 |
| | Real number part as long integer | 4888520323532708650 |
| | Imaginary number part as long integer | 4884430403359071803 |
| | Modulus of PC as long integer | 4889396296485818748 |
| | Theta of PC as long integer | 4602468698391823727 |
| | Cantor PF value | 45047750540491773913433549502792707777 |
| | Szudzik PF value | 23906196144089240399724999666785929979 |
| | Real number part as bit-string | 0100001111010111100000110101010101110001011100011101101100101010 |
| | Imaginary number part as bit-string | 0100001111001000111110111001001010110010010101011101111000111011 |
| | Modulus of polar coordinate as bit-string | 0100001111011010101000000000001101100101000100101110100010111110 |
| | Theta of polar coordinate as bit-string | 0011111111011111010000001110010101110001011111111110001010110111 |
| | Polar bit-interleaved value | 49679650227602418166657328428407275253 |
| | Cartesian bit-interleaved value | 63885745057879574985027357472836160973 |
| Polar Cantor value | Modulus of PC as long integer | 4889396296485818748 |
| | Theta of PC as long integer | 4602468698391823727 |
| | Modulus of PC | 7.67416362618991E18 |
| | Theta of PC | 0.4883359535588942 |
| | Real number part | 6.7771673222051697E18 |
| | Imaginary number part | 3.6003875414142126E18i |
| | *Error* | *6.671736816409259E-17* |
| Polar Szudzik value | Modulus of PC as long integer | 4889396296485818748 |
| | Theta of PC as long integer | 4602468698391823727 |
| | Modulus of PC | 7.67416362618991E18 |
| | Theta of PC | 0.4883359535588942 |
| | Real number part | 6.7771673222051697E18 |
| | Imaginary number part | 3.6003875414142126E18i |
| | *Error* | *6.671736816409259E-17* |
| Polar bit-interleaved value | Modulus of PC as bit-string | 0100001111011010101000000000001101100101000100101110100010111110 |
| | Theta of PC as bit-string | 0011111111011111010000001110010101110001011111111110001010110111 |
| | Modulus of PC | 7.67416362618991E18 |
| | Theta of PC | 0.4883359535588942 |
| | Real number part | 6.7771673222051697E18 |
| | Imaginary number part | 3.6003875414142126E18i |
| | *Error* | *6.671736816409259E-17* |
| Cartesian bit-interleaved value | Real number part as bit-string | 0100001111010111100000110101010101110001011100011101101100101010 |
| | Imaginary number part as bit-string | 0100001111001000111110111001001010110010010101011101111000111011 |
| | Real number part | 6.7771673222051697E18 |
| | Imaginary number part | 3.6003875414142131E18i |
| | *Error* | *0.0* |

Table 2
Representation errors of the proposed techniques

| PF | Maximum Error | Average Error |
|---|---|---|
| Polar Cantor | 1.04654252576E-15 | 2.04149108976E-16 |
| Polar Szudzik | 1.04654252576E-15 | 2.04149108976E-16 |
| Polar bit-interleaved | 1.04654252576E-15 | 2.04149108976E-16 |
| Cartesian bit-interleaved | 0.0 | 0.0 |

Based on the results shown in Table 2, the Cartesian bit-interleaved produces the best representation value by scoring 0.0 errors in both maximum and average representation error, as expected. It is also worth mentioning that all polar-based representation produces the same maximum and average error results, thus it is concluded that all polar-based representation carries the same discriminative power. The reason is that the conversion of Cartesian into polar coordinates may result in precision lost, and thus when converting back from polar into Cartesian coordinates, the original value will be slightly different.

To further corroborate the merit of Cartesian bit-interleaved as opposed to other polar-based techniques, thorough statistical validation using analysis of variance (ANOVA) should be performed. However, since the variance of Cartesian bit-interleaved is equals to zero ($\sigma = 0.0$) and the variances of polar-based techniques is equals to each other, the statistical validation is rendered impractical.

This study has successfully proposed four complex number representation techniques, with the Cartesian bit-interleaved is considered as the best proposed technique among other proposed techniques. This study also managed to achieve the objective of reducing the space required to retain the information. On the other hand, the validation of the proposed techniques using ML algorithms, most notably in the pattern recognition domain, will be conducted in the future works. This study believes that the findings have wide range of applicability, for example the representation of features extracted using Fast Fourier Transform, representation of frequency and amplitude for voice recognition, representation of complex-valued moments, such as Zernike, Chebyshev–Fourier, and orthogonal Fourier–Mellin moments, for 2D and 3D image analysis, etc.

## V. Conclusion and Future Works

Novel complex number representation techniques have been presented in this study, which are polar Cantor, polar Szudzik, polar bit-interleaved, and Cartesian bit-interleaved complex number representations. This paper also compared the merits of the proposed techniques. The experiments have shown that Cartesian bit-interleaved is the best representation technique for complex number, without any loss of information.

Hence, future works to incorporate Cartesian bit-interleaved in ML algorithms and toolkits, such as providing the implementation in WEKA, is planned. The proposed techniques and its Java source code are also publicly available in http://ftmk.utem.edu.my/cit/research/ats-drugs/tools.

## Acknowledgment

This work was supported by Universiti Teknikal Malaysia Melaka (UTeM) through the PJP High Impact Research Grant [PJP/2016/FTMK/HI3/S01473] and UTeM Postgraduate Fellowship (Zamalah) Scheme.